\documentclass[final,onefignum,onetabnum]{siamart171218}

\usepackage{cite}
\usepackage{graphicx}
\usepackage{amsmath,amsfonts}
\usepackage{epsfig}
\usepackage{multirow}
\newcommand{\bx}{{\bf x}}

\newcommand{\bd}{{\bf d}}
\newcommand{\bh}{{\bf h}}
\newcommand{\bb}{{\bf b}}

\newcommand{\be}{{\bf e}}

\usepackage{lipsum}
\usepackage{amsfonts}
\usepackage{graphicx}
\usepackage{epstopdf}
\usepackage{algorithmic}
\ifpdf
  \DeclareGraphicsExtensions{.eps,.pdf,.png,.jpg}
\else
  \DeclareGraphicsExtensions{.eps}
\fi

\newsiamremark{remark}{Remark}
\newsiamremark{hypothesis}{Hypothesis}
\crefname{hypothesis}{Hypothesis}{Hypotheses}
\newsiamthm{claim}{Claim}

\headers{$L_1/L_2$ Minimization in Single-View Tomography}{}

\title{Box-Constrained $L_1/L_2$ Minimization in Single-View Tomographic Reconstruction}

\author{Sean Breckling
          \thanks{Nevada National Security Site  
          (\email{brecklsr@nv.doe.gov}, \email{pillowja@nv.doe.gov}, \email{baldonbj@nv.doe.gov}).}
\and Malena I. Espa\~nol 
\thanks{Arizona State University
  (\email{mespanol@asu.edu}, 
  \email{vduribe@asu.edu}).}
\and Victoria Uribe\footnotemark[2]
\and Chrisitan Bobmara\thanks{Yale Univesity (\email{christian.bombara@yale.edu}).}
\and Jordan Pillow\footnotemark[1]
\and Brandon Baldonado\footnotemark[1]}

\usepackage{amsopn}

\makeatletter
\newcommand*{\addFileDependency}[1]{% argument=file name and extension
  \typeout{(#1)}% latexmk will find this if $recorder=0 (however, in that case, it will ignore #1 if it is a .aux or .pdf file etc and it exists! if it doesn't exist, it will appear in the list of dependents regardless)
  \@addtofilelist{#1}% if you want it to appear in \listfiles, not really necessary and latexmk doesn't use this
  \IfFileExists{#1}{}{\typeout{No file #1.}}% latexmk will find this message if #1 doesn't exist (yet)
}
\makeatother

\newcommand*{\myexternaldocument}[1]{%
    \externaldocument{#1}%
    \addFileDependency{#1.tex}%
    \addFileDependency{#1.aux}%
}
\myexternaldocument{ex_supplement}

\begin{document}

\maketitle

% REQUIRED
\begin{abstract}
 We present a note on the implementation and efficacy of a box-constrained $L_1/L_2$ regularization in numerical optimization approaches to performing tomographic reconstruction from a single projection view. The constrained $L_1/L_2$ minimization problem is constructed and solved using the Alternating Direction Method of Multipliers (ADMM). We include brief discussions on parameter selection and numerical convergence, as well as detailed numerical demonstrations against relevant alternative methods. In particular, we benchmark against a box-constrained TVmin and an unconstrained Filtered Backprojection in both cone and parallel beam (Abel) forward models. We consider both a fully synthetic benchmark, and reconstructions from X-ray radiographic image data. 
\end{abstract}

% REQUIRED

% REQUIRED

\section{Introduction}
\label{intro}
\begin{figure}
\begin{center}
Idealized Imaging System Diagram
\includegraphics[width=4in]{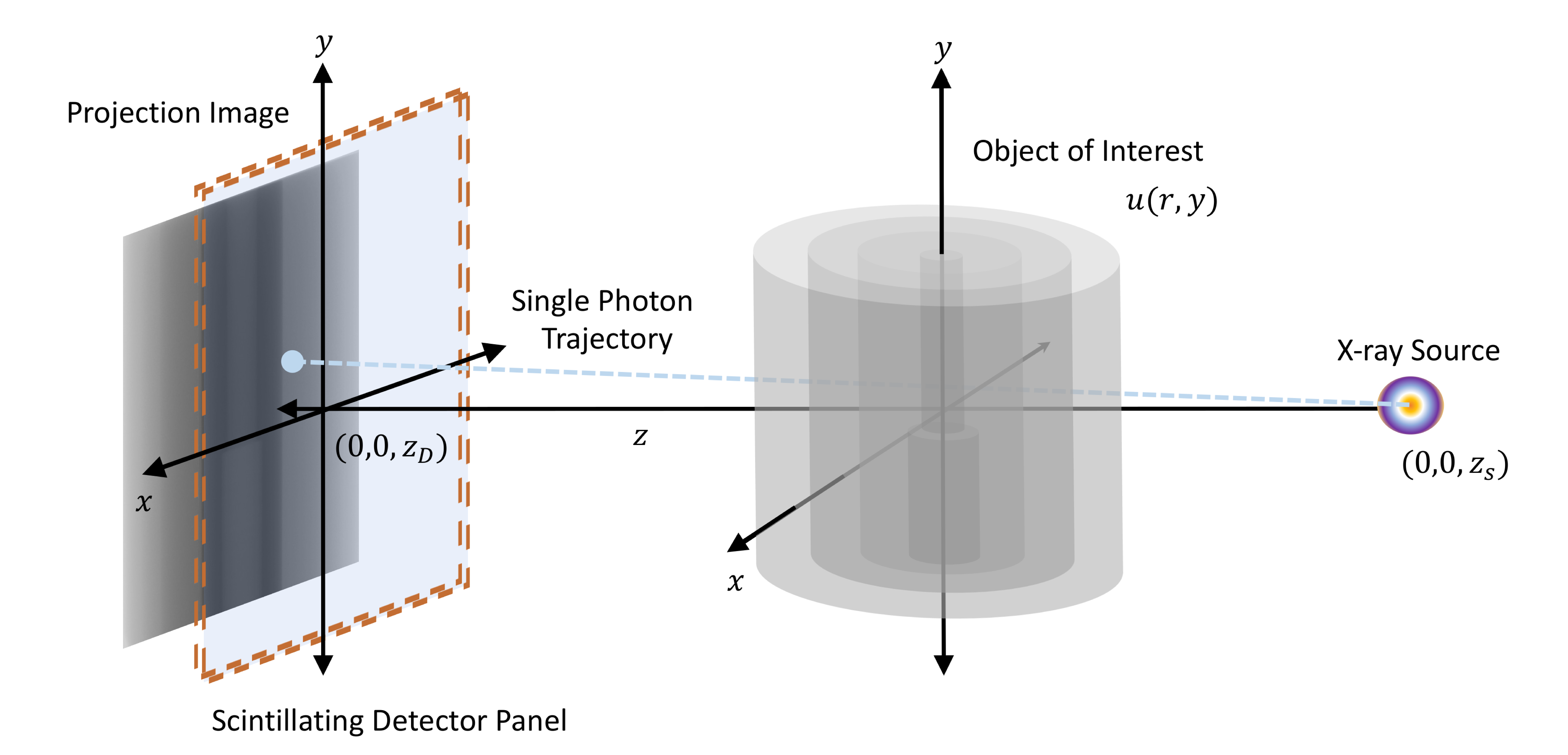}
 \end{center}
 \caption{\label{fig::2DDiag}A diagram depicting the spatial orientation of an X-ray imaging system in 3D. The dashed line represents a hypothetical photon path, $u(r,y)$ a target object axially-symmetric about the $y$-axis, $(z_s,0, 0)$ the location of an X-ray source, and $(z_D, x, y)$ the point along the photon path intersecting with the X-ray detector plane.}
\end{figure}

The reconstruction of three-dimensional information from a single projection image is a common task in quantitative X-ray radiography \cite{vest1974formation,hanson1985tomographic}, phase microscopy \cite{roberts2002refractive}, pyrometry \cite{wang1978temperature, dreyer2019improved}, and observational astronomy \cite{hansen1978optical, hasenberger2020aviator}. When the object of interest is known to present strong axial symmetry about an axis ($y$) perpendicular to the central beam propagation axis ($z$), the reconstruction problem becomes equivalent to solving the associated integral equation 
\begin{equation}\label{eqn::integral_eqn}
    d(x,y) = A\left\lbrace u(r,y)\right\rbrace = \int_{|a(x,y)|}^\infty \frac{r u(r,y)}{\sqrt{r^2 - a^2(x,y)}} dr,
\end{equation}
where $r = \sqrt{z^2 + x^2}$,  $d(x,y)$ is the resulting image recorded on the detector plane, $u(r,y)$ is the volumetric density of the target object of interest, and $a(x,y)$ denotes the minimum distance between a particular X-ray path and the axis of symmetry. A 3D mock-up is shown in Figure \ref{fig::2DDiag}.

Reconstructing volumetric information in 3D from a single projection view comes with a host of caveats. First, the operator $A$ in \eqref{eqn::integral_eqn} is unbounded. This means that discretizations of $A$ must be treated with care. Second, for an imaging scientist to make effective use of this formulation, the X-ray apparatus and target object need to strictly adhere to the symmetry assumption. Even deviations within otherwise practical uncertainty limits can result in reconstructions presenting severe error accumulation near the axis of symmetry. 

% 0.85/7 , -7/0.85, mag =  b=5*0.85/7, (x + 5) = -(7/0.85)^2 x

Several techniques have been developed to address the shortcomings of the classic least-squares, or filtered back-projection approaches to solving (\ref{eqn::integral_eqn}). Some examples include multi-point approximations of the inverse operator, often referred to as deconvolution methods \cite{dasch1992one}, Kalman filter methods \cite{Hansen:85}, and basis set expansion (BASEX) \cite{dribinski2002reconstruction}. In more recent years, so-called regularized methods have been shown to be a reliable, robust alternative \cite{asaki2005abel,Asaki2006,chan2015high,howard2016bayesian}. Unconstrained regularization methods solve \eqref{eqn::integral_eqn}  by minimizing 
\begin{equation}
\min_{u} \ \mathcal{R}(u) + \frac{\lambda}{2}||Au - d||_2^2, 
\end{equation}
where $\mathcal{R}(u)$ is a non-negative regularizing penalty, and $\lambda > 0$ is a real-valued weighting parameter. Among the most well-known of these methods is Total-Variation Minimization (TVmin) \cite{rudin1992nonlinear}, which select
\begin{equation} \label{eqn::tvmin}
    \mathcal{R}(u) = || \nabla u ||_1. 
\end{equation}

In image restoration problems, choosing (\ref{eqn::tvmin}) over, say, the Tikhonov method  $\mathcal{R}(u) = ||\nabla u||_2$, is shown to penalize numerical optimization algorithms away from results presenting large amounts of non-sparsely distributed variation without over-smoothing prominent edges. Another benefit is the relative ease in which a full optimization procedure can be directly implemented \cite{chambolle2011first, goldstein2009split, eckstein1992douglas}. In literature, the primary shortcoming mentioned of (\ref{eqn::tvmin}) are “stair-case” artifacts, which result when the regularization penalty plays an outsized role during minimization. 

Given the wide, general applicability of the TVmin approach to sparse signal recovery problems, numerous alternative formulations of $\mathcal{R}(u)$ have been developed. The focus of this paper are recent results from \cite{wangnagy2021} demonstrating promising results from the  $L_1/L_2$ regularizer 
\begin{equation}\label{eqn::l1_l2}
    \mathcal{R}(u) = \frac{|| \nabla u ||_1}{|| \nabla u ||_2}
\end{equation}
\noindent in limited-angle tomography problems. Given that solving (\ref{eqn::integral_eqn}) is the extreme case of limited-angle tomography, we are keenly motivated to test if similar improvements can be seen in this modality.

This study implements a numerical optimization approach to minimizing the constrained model
\begin{equation}\label{eqn:main_prob}
  \min_{u}  \  \frac{|| \nabla u ||_1}{|| \nabla u ||_2} +  \frac{\lambda}{2}||Au - d||_2^2 \ \text{ such that } \  \alpha \leq u \leq \beta,
\end{equation}
using the alternating direction method of multipliers (ADMM) \cite{glowinski1975approximation} technique implemented in \cite{wangnagy2021}, but using parallel beam (Abel), and cone-beam formulation of the axilly symmetric single-projection operation (\ref{eqn::integral_eqn}). In particular, we consider constraints like $u \geq 0$, given that the reconstructed objects represent non-negative physical quantities, such as the volumetric density of the target scene as a function of space. We also include a brief discussion on this method's numerical convergence in the single-view tomography settings. 

We provide two numerical demonstrations of an ADMM implementation of \eqref{eqn:main_prob} against an ADMM implementation of a similarly-constrained $TVmin$ model, as well as an unconstrained filtered back projection technique. These tests are conducted in two specific settings: the first compares performance in the parallel-beam modality against fully manufactured data, and the second in a cone-beam modality associated with radiography collected at the Scintillator Evaluation and Assessment Laboratory (SEALab) at the Los Alamos, New Mexico operations office of the Nevada National Security Site (NNSS). In performing these numerical experiments we found that, when properly parameterized, $L_1/L_2$ regularization is effective at reducing error near the singularity at $r=0$ in \eqref{eqn::integral_eqn}, as well as reducing artifacts common to over-regularization in TVmin or Tikhonov-regularized models. 

The remainder of this paper proceeds as follows. In Section 2, we discuss the mathematical settings of the reconstruction problem and introduce the main algorithm. A brief discussion of the convergence behavior and parameter selection is also included. In Section 3, we provide numerical demonstrations using the synthetic and genuine radiographic image data. Concluding remarks and discussion of future work are given in Section 4.

\section{Preliminaries and Numerical Methods}
We fix the X-ray source of our imaging device at $(x,y,z) = (0,0,z_S)$ and the flat, scintillating X-ray detector plate on the plane $(x,y,z_D)$. The 3D region bounding the object of interest is compactly supported between the source and detector, within the rectangular prism $\Omega = \mathcal{I}_x \times \mathcal{I}_y \times \mathcal{I}_z$, where $\mathcal{I}_{(\cdot)}$ denotes a closed interval containing the origin.  Since we assume, without loss of generality, that the principal X-ray beam axis is $z$ and that the object of interest is centered at the origin and axially symmetric about the $y$ axis, $\Omega$ can be reduced to the cylindrical coordinate system $\Phi = \mathcal{I}_r \times \mathcal{I}_y$ where $r^2 = x^2 + z^2$. 

Let $A : \Phi \rightarrow \Psi$ be the Abel operator, where $\Psi$ is on the surface of the detector plate $\mathcal{I}_x \times \mathcal{I}_y$ located at $z = z_D$. Let $\Phi_h$ and $\Psi_h$ denote uniform discretizations of the detector plate $\mathbb{R}^{n \times m}.$  We index discrete digital images in the usual way, where $d_{i,j} \in \mathbb{R}$ is shorthand for the $((i-1)m + j)$-th term in $d$. We also adopt the usual bold notation for vector-valued terms, e.g., $\bd \in \Psi_h \times \Psi_h$. We discuss discretizations to the integral operator \eqref{eqn::integral_eqn} in more detail below in Section \ref{sec::abel}.

Let $X$ denote either $\Phi_h$ or $\Psi_h$. We employ the usual inner product and $L_2$ norm notation throughout, i.e., 
\begin{equation}
(f,g)_{X} =  \sum_{i,j = 1}^{n,m} f_{i,j} g_{i,j} \mbox{ and } ||f||_2^2 = (f,f)_{X}. \nonumber
\end{equation}
We similarly denote the $L_1$ norm over $X$ in the usual way
\begin{equation}
    ||f||_1 = \sum_{i,j = 1}^{n,m} |f_{i,j}|. \nonumber
\end{equation}

We write the discrete, forward-differencing gradient operator as 
\begin{equation}
    \nabla f = \left( \nabla_r f, \nabla_y f \right). \nonumber 
\end{equation}
 Reflecting boundary conditions are enforced, which play a modest role in mitigating error accumulation near the $r=0$ axis. We define our discrete Laplace operator in terms of our discrete gradient such that $\triangle = -\nabla^T \nabla.$ We denote the anisotropic TV norm  $||f||_{TV(\Phi_h)} = ||\nabla f||_{L_1(\Phi_h \times \Phi_h)}$ as simply $||\nabla f||_1$, where
\begin{equation}
    ||\nabla f||_1 = \sum_{i,j = 1}^{n,m} \bigg( |\nabla_r f_{i,j} | + |\nabla_y f_{i,j}|\bigg) = ||\nabla_r f ||_1 + ||\nabla_y f||_1. \nonumber
\end{equation}
Lastly, we adopt the following notation to denote the usual shrink function
$${\bf shrink}(\bx,\mu) = \text{sign}(\bx)\max\lbrace |\bx| - \mu, 0\rbrace.$$
\subsection{Discretizing The Integral Transform} \label{sec::abel} 
Again consider the diagram in Figure \ref{fig::2DDiag}. We denote the distance from the X-ray source to the object center as $|z_s|$, and the distance from the object center to the detector plane as $|z_D|$. The apparent magnification of the imaging system is given by 
\begin{equation*}
    \xi = \frac{|z_D| + |z_S|}{|z_S|}.
\end{equation*}
Pixel spacing at the focal plane $\Psi_h$ is scaled proportionately with the magnification of the imaging system such that $\Delta r= \Delta x / \xi$. The same scaling occurs along the $y$ axis, and thus we make the distinction $\bar{y} = y/\xi$.

In $\mathbb{R}^3$, given a particular pixel $(i,j)$ on at the detector $\Psi_h$, physically located at $(x_i,y_j,z_D)$, the ideal X-ray photon path is the line parameterized for $t\in[0,1]$:
\begin{equation}\label{eqn::line_seg}
    L(i,j)(t) = (0,0,z_s) + t(x_i, y_j, z_D - z_S).
\end{equation}
Next, let $(i',j')$ denote pixels in $\Phi_h$, which each represent an annulus of the form $\mathcal{R}_{i',j'} = [r_{i'}, r_{i'+1}] \times [\bar{y}_{j'},\bar{y}_{j'+1}]$. 

If we assume that our discretized reconstructions are of the form
\begin{equation*}
    u(x,y) = \sum_{\forall i,j} u_{i,j} \chi_{\mathcal{R}_{i,j}},
\end{equation*} where $\chi_{\mathcal{R}_{i,j}}$ denotes the characteristic function over the $(i,j)^{th}$ annulus, then it follows from a direct analytic computation of $A\lbrace{u(x,y)}$ that the discretized projection operator $A : \Phi_h \rightarrow \Psi_h$ can be defined by assigning each particular coefficient $A_{(i,j), (i',j')}$ to be the length $L_{i,j} \cap \mathcal{R}_{i',j'}.$ In these general settings, this method results in a cone-beam formulation of the onion-layer type discretizations in wide use for parallel-beam modalities \cite{dasch1992one, pretzier1992}. 

\subsection{An ADMM for the Box-Constrained Problem}
Given that we have made no procedural adjustments to the ADMM developed in \cite{wangnagy2021}, we provide only the details relevant to the implementation for solving \eqref{eqn:main_prob}. While selecting this particular implementation necessitates three iteration levels (outer and inner ADMM loops, as well as an innermost Krylov solve), single-view forward models are rarely as computationally expensive as those for traditional, limited-angle, or even few-view computed tomography. Therefore, the added cost is less of a practical consideration.

\begin{algorithm}
\caption{The Box-Constrained $L_1/L_2$ Minimization}\label{alg:box}
{\bf Input:} The projection operator $A$, observed data $d$, and reconstruction bounds $[\alpha,\beta]$.  \\
{\bf Parameters:} $\lambda, \rho_1, \rho_2, \rho_3, \varepsilon \in \mathbb{R}^+,$ and  $k_{max}, j_{max} \in \mathbb{N}.$ \\
{\bf Initialize:} $\bh, \bb_1, \bb_2, \bd = {\bf 0}$, and $e,k,j = 0.$ \\
{\bf Initialize:} $M = \lambda A^T A - (\rho_1 + \rho_2)\triangle + \rho_3 I.$
\begin{algorithmic}
\WHILE{ $k < k_{max} $ \ or $||u^{(k)}-u^{(k-1)}||_2 / ||u^{(k)}||_2 > \varepsilon $}

\WHILE{$j < j_{max}$ \ or $||u_j-u_{j-1}||_2 / ||u_j||_2 > \varepsilon$}
\STATE $u_{j+1} = M^{-1}(\lambda A^T d + \rho_1\nabla^T(\bd_j - (\bb_1)_j) +  \rho_2\nabla^T(\bh^{(k)}-\bb_2^{(k)}) + \rho_3(v^{(k)}-e^{(k)}))$
\STATE $\bd_{j+1} = {\bf shrink}\left(\nabla u_{j+1}+(\bb_1)_j, \frac{1}{\rho_1 || \bh^{(k)}||_2} \right)$
\STATE $v_{j+1} = \min\lbrace \max \lbrace u_{j+1}+e_j, \alpha \rbrace, \beta \rbrace$
\STATE $(\bb_1)_{j+1} = (\bb_1)_j + \nabla u_{j+1} - \bd_{j+1}$
\STATE $e_{j+1} = e_j +u_{j+1} - v_{j+1}$
\STATE $j = j+1$
\ENDWHILE
\STATE $u^{(k+1)} = u_j$
\STATE $\bh^{(k+1)} = \begin{cases} \tau^{(k)}\left( \nabla u^{(k+1)} + \bb_2^{(k)}\right) \ &\text{ if }  \nabla u^{(k+1)} + \bb_2^{(k)} \neq {\bf 0} \\ \be^{(k)} \ &\text{ otherwise.} \end{cases}$
\STATE $\bb_2^{(k+1)} = \bb_2^{(k)} + \nabla u^{(k+1)} - \bh^{(k+1)}$
\STATE $ k = k+1$
\STATE $j = 0$
\ENDWHILE
\end{algorithmic}
{\bf return} $u^* = u^{(k)}$
\end{algorithm}

For $\bh$-update, we require either 
$$\tau^{(k)} = \frac{1}{3}\left( C^{(k)} + 1 + \frac{1}{C^{(k)}} \right),$$  
where
$$C^{(k)} = \sqrt[3]{\frac{27D^{(k)} + 2 + \sqrt{(27 D^{(k)}+2 )^2-4}}{2}},$$ 
and
$$D^{(k)} = (\rho_2 ||\nabla u^{(k+1)} + \bb_2^{(k+1)} ||_2^3)^{-1} ||\nabla u^{(k+1)}||_1,$$
or a uniformly-distributed random array $\be^(k)$ such that $$||\be^{(k)}||_2^3 = \rho_2^{-1} ||\nabla u^{(k+1)}||_1.$$

What follows guarantees numerical convergence of Algorithm \ref{alg:box}.
 \begin{theorem}[Convergence of Algorithm \ref{alg:box}] Given the axially-symmetric forward operator $A$ defined in Section \ref{sec::abel}, if $Null(\nabla) \cap Null(A) = \lbrace {\bf 0} \rbrace$ and the sequence $\lbrace || \bh^{(k)} ||_2 \rbrace$ permits a positive, uniform lower bound, then for a sufficiently large $\rho_2$, the results of the outermost loop sequence $\lbrace u^{(k)},\bh^{(k)} \rbrace$ always admits a subsequence convergent to a critical point of \eqref{eqn:main_prob}.
 \end{theorem}
 
 Given the first assumption represents the principal departure from the work in \cite{wangnagy2021}, we refer the reader to Theorem 4.6 therein. It is indeed verifiable that the single-view forward-operator $A$ satisfies the presumed conditions directly, regardless of the boundary conditions assigned to $\nabla$, since that $Null(A) = \lbrace{\bf 0}\rbrace$. 
 
 \begin{remark}[Parameter Selection] \label{tuning}
Selecting appropriate values for $\lambda, \rho_1, \rho_2$, and $\rho_3,$ is not often a straightforward task. If we first consider the inner update for $u$, we notice that one could scale $\lambda$ according to $||A^TA||$, whereas $\rho_1$ and $\rho_2$ scale with $||\triangle||,$ and $\rho_3$ with $I.$ In our numerical tests that follow, we found that selecting $\lambda \approx \mathcal{O}(1) ||A^TA||^{-1}$, $\rho_1$ and $\rho_2 \approx \mathcal{O}(1) \Delta r^2$ and $\rho_3 \approx \mathcal{O}(1)$ was a serviceable starting point for tuning. A more-general guidance can be found in \cite{tao2022minimization}. 
\end{remark}

\section{Numerical Demonstrations}
This section exhibits two implementations of Algorithm \ref{alg:box}. We compare those implementations against the default GRIDREC Fourier-based solver \cite{dowd1999developments} in the \texttt{tomopy} package \cite{gursoy2014tomopy}, and an ADMM implementation of a min-max constrained TVmin \cite{chan2013constrained} (Algorithm \ref{alg:boxtv}). The first test considers synthetic data, wherein we generate known ground-truth values for both the target object and the projection. The second demonstration performs a similar comparison, but utilizes X-ray radiography collected at SEALab.

\begin{algorithm}[H]
\caption{The Box-Constrained TV Minimization}\label{alg:boxtv}
{\bf Input:} The projection operator $A$, observed data $d$, and reconstruction bounds $[\alpha,\beta]$.  \\
{\bf Parameters:} $\lambda, \rho_1, \rho_2, \varepsilon \in \mathbb{R}^+,$ and  $j_{max} \in \mathbb{N}.$ \\
{\bf Initialize:} $\bh, \bb = {\bf 0}$, and $e,j = 0.$ \\
{\bf Initialize:} $M = \lambda A^T A - \rho_1\triangle + \rho_2 I.$
\begin{algorithmic}

\WHILE{$j < j_{max}$ \ or $||u_j-u_{j-1}||_2 / ||u_j||_2 > \varepsilon$}
\STATE $u_{j+1} = M^{-1}(\lambda A^T d +  \rho_1\nabla^T(\bh_j-\bb_j) + \rho_2(v_j-e_j))$
\STATE $\bh_{j+1} = {\bf shrink}\left(\nabla u_{j+1}+\bb_j, \rho_1^{-1} \right)$
\STATE $\bb_{j+1} = \bb_j + \nabla u_{j+1} - \bh_{j+1}$
\STATE $v_{j+1} = \min\lbrace \max \lbrace u_{j+1}+e_j, \alpha \rbrace, \beta \rbrace$
\STATE $e_{j+1} = e_j +u_{j+1} - v_{j+1}$
\STATE $j = j+1$
\ENDWHILE
\end{algorithmic}
{\bf return} $u^* = u^{(k)}$
\end{algorithm}

\subsection{Comparisons with Synthetic Targets and Projections} 
\begin{figure}
    \centering
    \includegraphics[height=2in]{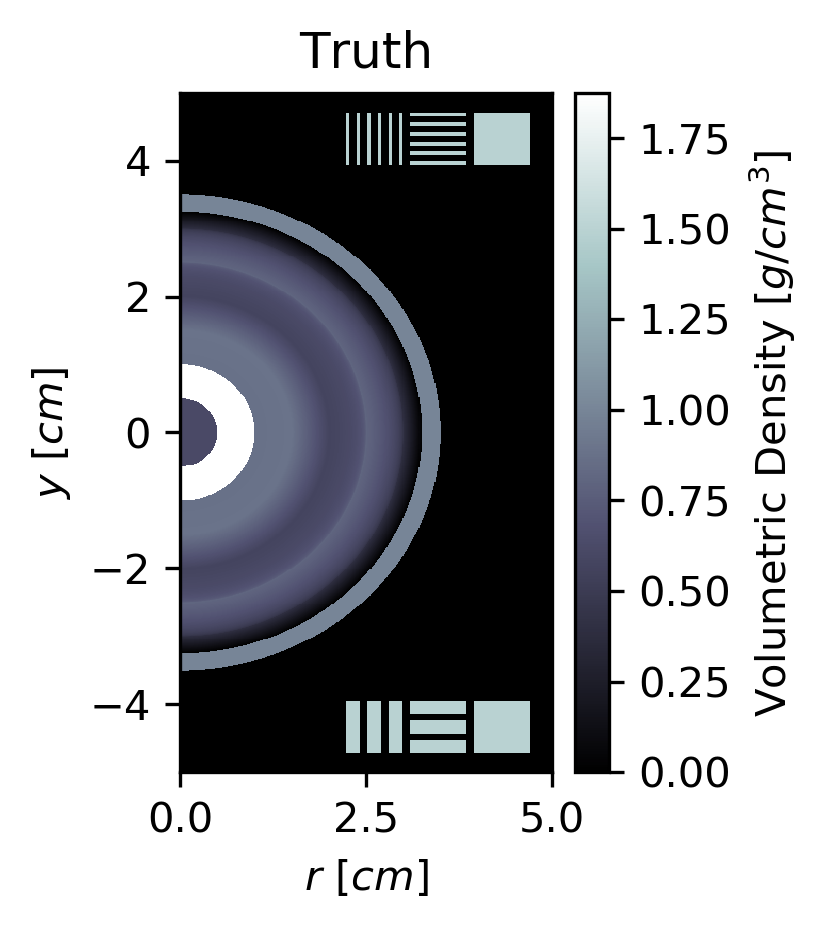}\includegraphics[height=2in]{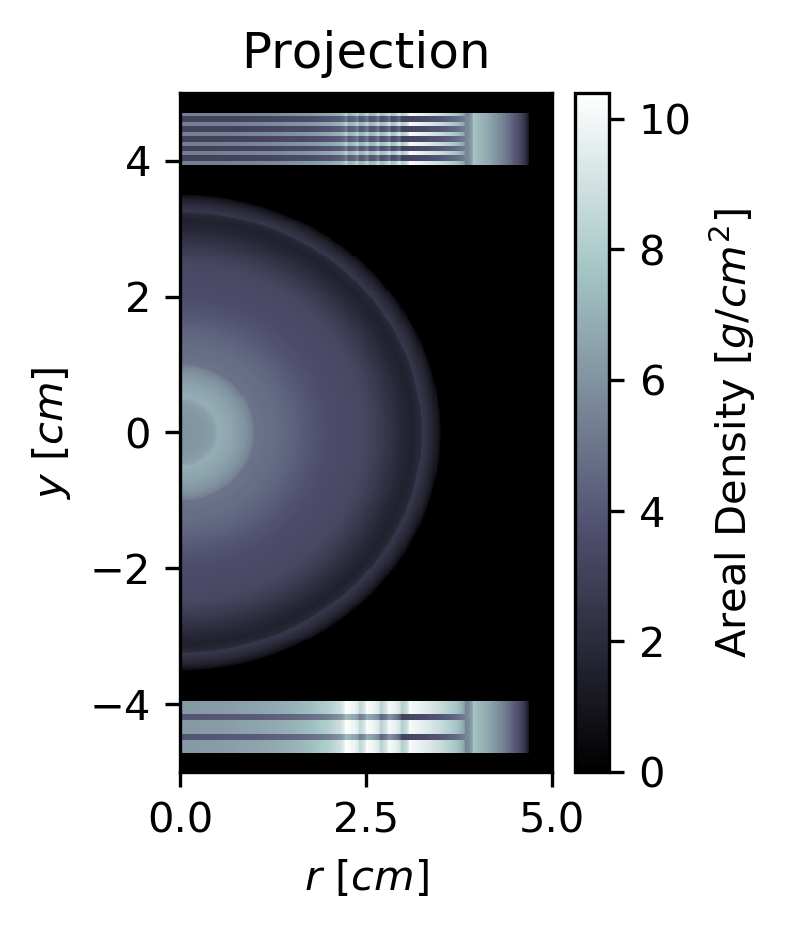} \caption{Synthetic data consisting of linear combinations of spherically-symmetric functions centered at the origin, and constant-valued annular fiducial objects at the top and bottom right corners (Left). Given that this synthetic data is constructed with functions selected from Table \ref{tab:test_functions}, we also have a ground truth of the forward projection (Right).
    \label{fig:my_label}}
\end{figure}

Our first test is conducted using synthetic radiography, which can be seen in Figure \ref{fig:my_label}. The inner, spherically-symmetric functions in Figure \ref{fig:my_label} are constructed to be invariant along concentric spheres centered at the origin of the $(r,y)$ half-plane. In cylindrical coordinates, we use the mapping
$$\delta(y, \nu) = \nu \sin \left( \arccos(y/\nu)\right) $$ in conjunction with the functions in Table \ref{tab:test_functions},
where $\nu>0$ is a prescribed spherical radius. 

{\renewcommand{\arraystretch}{2}%

\begin{table}[htbp]
{\footnotesize
    \caption{Linear combinations of the following functions can be used to generate synthetic axially-symmetric data. }
    \centering
    \begin{tabular}{|c|c|c|}
       \hline 
        $k$ & $f_k(r;\delta), \ r\leq \delta $             & $A\lbrace f_k\rbrace(x; \delta), \ x\leq \delta$ \\
         \hline                               
        1 & $ 1 $   & $2\sqrt{\delta^2 - x^2} $ \\
        2 & $ \sqrt{\delta^2 - r^2} $ & $\frac{\pi}{2}(\delta^2 - x^2) $ \\
        3 & $ (\delta^2 - r^2)^{3/2}$& $\frac{3\pi}{8}(\delta^2 - x^2)^2 $\\
        \hline 
         
    \end{tabular}

    \label{tab:test_functions}}
\end{table}}

Our synthetic ground truth function $u^*$ is constructed as
\begin{multline}\label{eqn:trueu}
    u^*(r,y) = \sum_{i=1}^4 \alpha_i f_1\left(r;\delta\left(y;\nu_i\right)\right) + \sum_{i=5}^7 \alpha_i f_2\left(r;\delta\left(y;\nu_i\right)\right) \\ + \sum_{i=8}^{12} \alpha_i f_3\left(r;\delta\left(y;\nu_i\right)\right) +  \sum_{i=1}^{34} \tfrac{3}{2} \chi_i(r,y),
\end{multline}
where 
\begin{align*}
    \alpha & = \left\lbrace  0.50, 1.00, 3.25, 3.50, 2.50, 3.00, 1.50, 2.00, 3.00, 3.25 \right\rbrace, \\
    \nu & = \left\lbrace -1.25, 1.00, -1.00, 1.00, -0.25, 0.25, -0.20, 0.20, -0.25, 0.25 \right\rbrace.
\end{align*}
Each $\chi_i$ corresponds to a particular rectangular characteristic annulus. These  features are generated with linear combinations of $f_1,$ varying fixed choices of $\delta$ to correspond to selected radii.  All rectangles have one side with length $l=0.75$ with an opposing side of length $l, 3l/11,$ or $l/11,$ oriented as displayed. The synthetic projections are constructed similarly, substituting $f_k$ with the corresponding $A\lbrace f_k \rbrace$ from Table \ref{tab:test_functions}, or by using \texttt{exact\_sinogram}  \cite{Dessole2023}.

Given our application, we assign metric length and volumetric density values in $cm$ and $g/cm^3,$ respectively. The reconstruction and projection domains $\phi_h$ and $\psi_h$ are both $[-5,5] \times [0,5]$ spaces discretized uniformly on an $N = 700\times 350$ pixel grid. As a result, $\Delta x = \Delta y = \Delta r = (1/70) cm.$

We constructed this function to exhibit challenging features common to the single-view tomography modality, to include noise at two levels: $\sigma = 0.25\%$ and $2.5\%$ of $||A(u^*)||_\infty$. Given that regularized methods often struggle to accurately reconstruct non piece-wise constant features, our central target displays a high degree of variability. Additionally, the reconstructions of the rectangular annuli can be used to identify losses in contrast or resolution. 

Since ground truth $u^*$ is known, the comparison below considers both the root mean squared error (RMSE) 
$$ RMSE(u,u') := \frac{||u-u'||}{N},$$ and the structural similarity (SSIM) 
 $$ SSIM(u,u^*) := \frac{1}{M} \sum_{i=1}^M \text{ssim}(u_i,u_i^*),$$
with $M$ denoting the number of discrete sub-image blocks considered. For this particular case, we select an exact tiling of sub-image squares of size $10\times10$ such that for two corresponding blocks $v, w$ $$\text{ssim}(v,w) := \frac{(2\mu_v \mu_w + c_1)(2\sigma_{vw} + c_2)}{(\mu_v^2 + \mu_w^2 + c_1)(\sigma_v^2 + \sigma_w^2 + c_2)}.$$

Parameter selection for Algorithms \ref{alg:box} and \ref{alg:boxtv} was accomplished by tuning as discussed in Remark \ref{tuning}. For Algorithm \ref{alg:box}, we noticed satisfactory results with $\lambda = 0.99||A^T A||^{-1}$ , $\rho_1 = \rho_2 = 5\times10^{-3} \Delta r^2$, and $\rho_3 = 1.0$. We set the maximum iterations to be $k_{max}=30$ and $j_{max}=5$, with a tolerance threshold of $\epsilon = \times 10^{-7}$. We select, effectively, the same choices for Algorithm \ref{alg:boxtv} with $\lambda = 0.99||A^TA||^{-1}$, $\rho_1 = 10^{-2} \Delta r^2$, and $\rho_2 = 1$. The maximal iteration count was fixed to be $j_{max} = 150,$ with a tolerance $\epsilon = 10^{-7}.$ The inner $u-$update step for both Algorithms \ref{alg:box} and \ref{alg:boxtv} is computed using the non-preconditioned Congjugate Gradient (CG) method from \texttt{scipy.sparse.linalg} \cite{2020SciPy-NMeth}. The maximum CG iterations were 1000, with the tolerance threshold set to $10^{-7}.$ We perform filtered backprojection by constructing a uniform $180^\circ$ sinogram with $\Delta \theta = 0.25^\circ$ from the single projection view.

\begin{figure}
    \centering
    \includegraphics[width=4.5in]{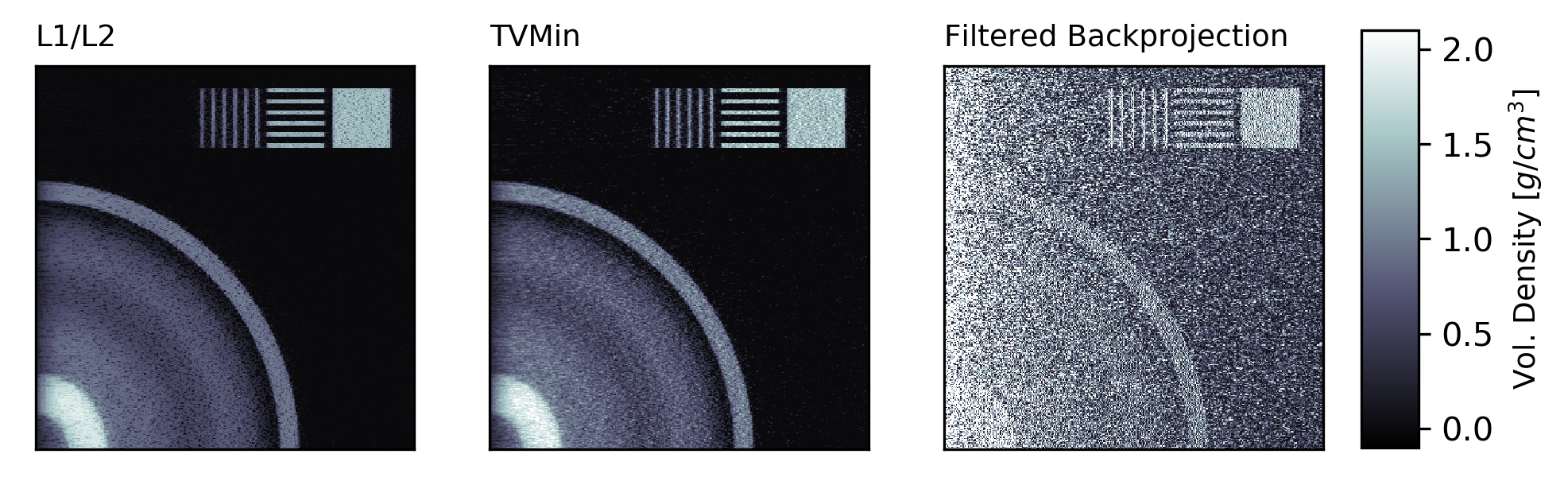} \caption{Three candidate solutions to the synthetic single-view reconstruction problem, where the i.i.d. additive noise parameter is selected as $\sigma = 0.025||A(u^*)||_\infty$. From left to right we have the results from Algorithm \ref{alg:box}, Algorithm \ref{alg:boxtv} (Center), and Filtered Backprojection (Right). Only the upper half-plane is shown to better highlight detail.
    \label{fig:synthrecon}}
\end{figure}

\begin{figure}
    \centering
    \includegraphics[height=3in]{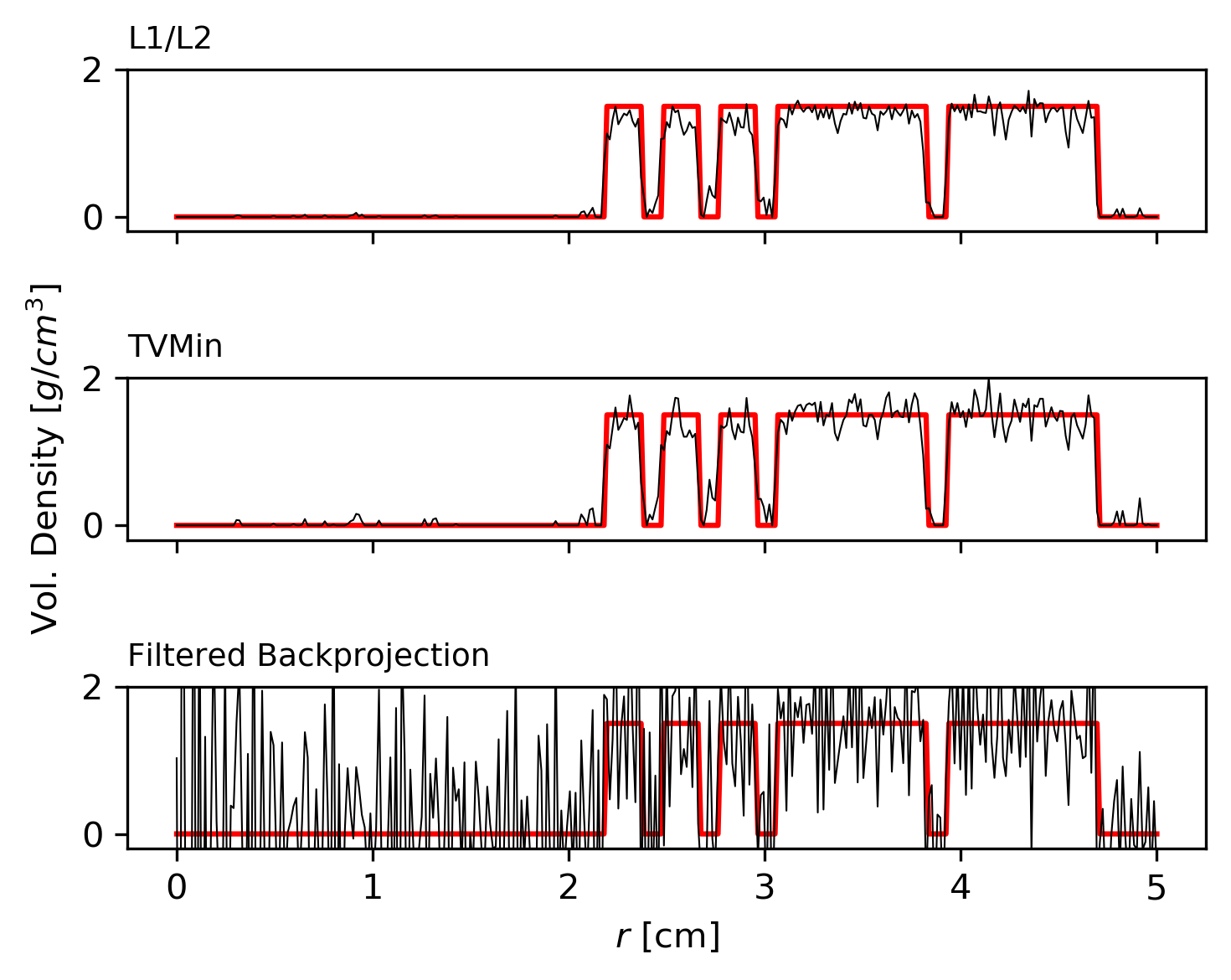} \caption{We compare the three candidate solutions to the synthetic single-view reconstruction problem along the line $y=-4.1.$ Similar to Figure \ref{fig:synthrecon}, the i.i.d. additive noise parameter is selected as $\sigma=2.5\times10^{-2}||A(u^*)||_\infty$. From top to bottom we have Algorithm \ref{alg:box} (Top), Algorithm \ref{alg:boxtv} (Middle), and Filtered Backprojection (Bottom).
    \label{fig:synthlineout}}
\end{figure}

In Figure \ref{fig:synthrecon} we show the upper half-planes of the results from all three methods in consideration. Figure \ref{fig:synthlineout} displays horizontal segments of the results through the rectangular features along the line  $y = -4.1$. In both figures, it is evident that regularization plays an enormous role in managing noise amplification. Further, we see that the $L_1/L_2$ method results in slightly-improved local accuracy; particularly in the empty regions. This is evidenced by the RMSE and SSIM scores in Table \ref{tab:rmsessim}.

\begin{table}[htbp]
{\footnotesize
  \caption{Comparing the RMSE and SSIM of reconstructions from Algorithms \ref{alg:box}, \ref{alg:boxtv}, and filtered backprojection.}
  \label{tab:rmsessim}
\begin{center}
\begin{tabular}{crccc}
\multicolumn{1}{l}{}                                & \multicolumn{1}{l}{}      & FBP                       & TVMin                    & $L_1/L_2$                \\ \hline
\multicolumn{1}{|c|}{\multirow{2}{*}{$\sigma=2.5 \times 10^{-3}||A(u^*)||_\infty$}} & \multicolumn{1}{r|}{RMSE} & \multicolumn{1}{c|}{$5.01\times 10^{-4}$}  & \multicolumn{1}{c|}{$2.99\times10^{-4}$}  & \multicolumn{1}{c|}{$2.90\times10^{-4}$}   \\ \cline{2-5} 
\multicolumn{1}{|c|}{}                              & \multicolumn{1}{r|}{SSIM} & \multicolumn{1}{c|}{$3.71\times 10^{-1}$}  & \multicolumn{1}{c|}{$8.84\times 10^{-1}$} & \multicolumn{1}{c|}{$9.51\times 10^{-1}$} \\ \hline
\multicolumn{1}{|c|}{\multirow{2}{*}{$\sigma=2.5 \times 10^{-2}||A(u^*)||_\infty$}} & \multicolumn{1}{r|}{RMSE} & \multicolumn{1}{c|}{$4.26\times 10^{-3}$} & \multicolumn{1}{c|}{$3.71\times10^{-4}$}  & \multicolumn{1}{c|}{$3.67\times10^{-4}$}   \\ \cline{2-5} 
\multicolumn{1}{|c|}{}                              & \multicolumn{1}{r|}{SSIM} & \multicolumn{1}{c|}{$3.69 \times 10^{-2}$}  & \multicolumn{1}{c|}{$4.91\times 10^{-1}$} & \multicolumn{1}{c|}{$6.25\times 10^{-1}$} \\ \hline
\end{tabular}
\end{center}
}
\end{table}

\subsection{Radiographic Image Data from SEALab}
The data presented below resulted from X-ray images of a cylindrical Aluminum calibration object machined to a $4cm$ outer diameter, with a symmetric cavitation machined away at various radii. We denote the X-ray image data as both $d(x,y)$ when referring to cartesian space, and $d_{i,j}$ when discussing discretizations. 

\begin{figure}
    \centering
    \includegraphics[height=3in]{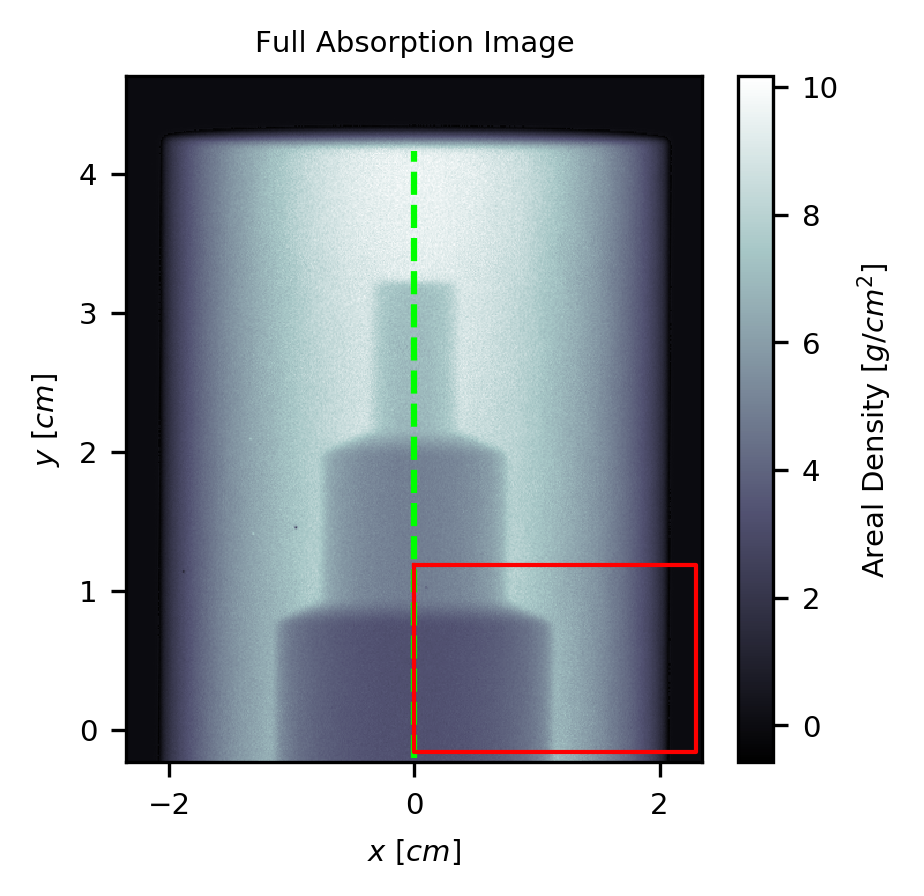}
    \caption{Areal density of Aluminum fiducial target object, resulting from applying a transfer curve to a flat-field normalized X-ray radiograph. The axis of symmetry is highlighted in green. We identify an ROI bound by the red box to better-highlight details.}
    \label{fig:projSealab}
\end{figure}

The object was positioned such that the shortest optical path between the source and detector intersected with the widest cavitation, which is $2cm$. We assign the intersection to be $(x,y) = (0,0)$ in $\Psi_h$, which is positioned toward the bottom in Figure \ref{fig:projSealab}. 

Selecting a threshold of $0.5g/cm^2$ such that $\Omega_1 = \lbrace d_{i,j}< 0.5\rbrace$ and $\Omega_2 = \lbrace d_{i,j} \geq 0.5 \rbrace$, the Contrast to Noise ratio (CNR) was measured as
\begin{equation*}
    CNR(d_{i,j}) = \frac{\left|\mu\left(d\left(\Omega_1\right) \right) - \mu\left(d\left(\Omega_2\right)\right) \right|}{\left|\sigma\left(d\left(\Omega_1\right) \right) - \sigma\left(d\left(\Omega_2\right)\right) \right|} = 7.564.
\end{equation*}
The calibration object was positioned in the field of view such that the vertical axis of symmetry in $\Phi_h$ corresponded to the $y$ axis. The distance from the origin $(0,0,0)$ to the distance to the detector plane measured to be $70.3cm$, and the corresponding distance to the X-ray source was measured to be $59.2cm$. The result is an apparent magnification of $1.1875.$ These distances were used in the construction of the cone-beam forward operator.
We parameterize Algorithms \ref{alg:box} and \ref{alg:boxtv} and the filtered backprojection identically to what was done in the previous subsection.

\begin{figure}
    \centering
    \includegraphics[width=5in]{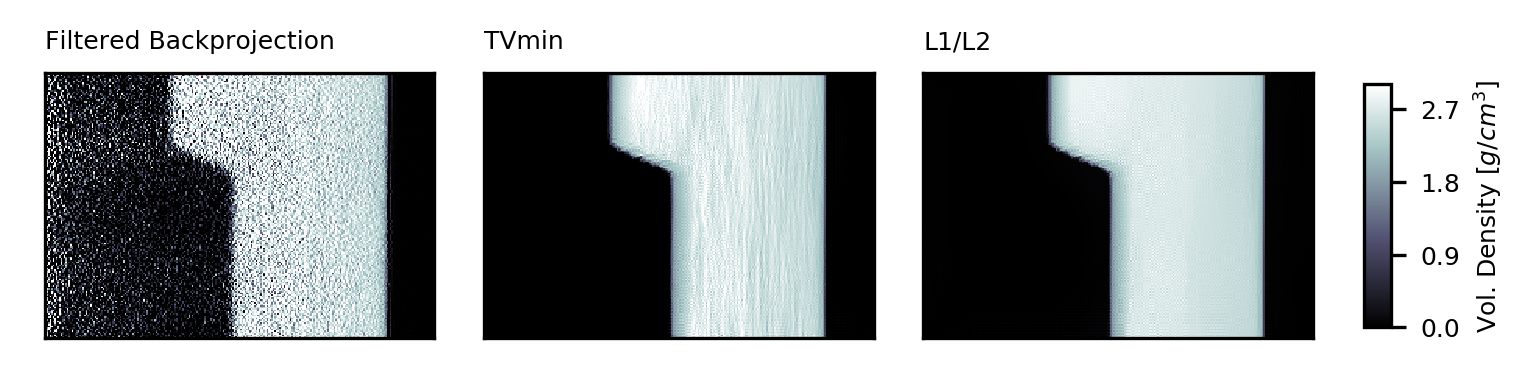}
    \caption{Volumetric density of Aluminum fiducial target object reconstructed from X-ray radiography using filtered backprojection (Left), TVmin (center) and $L_1/L_2$ (Right). All three images were reconstructed from the ROI highlighted in Figure \ref{fig:projSealab}.}
    \label{fig:compare3}. 
\end{figure}

\begin{figure}
    \centering
    \includegraphics[height=3in]{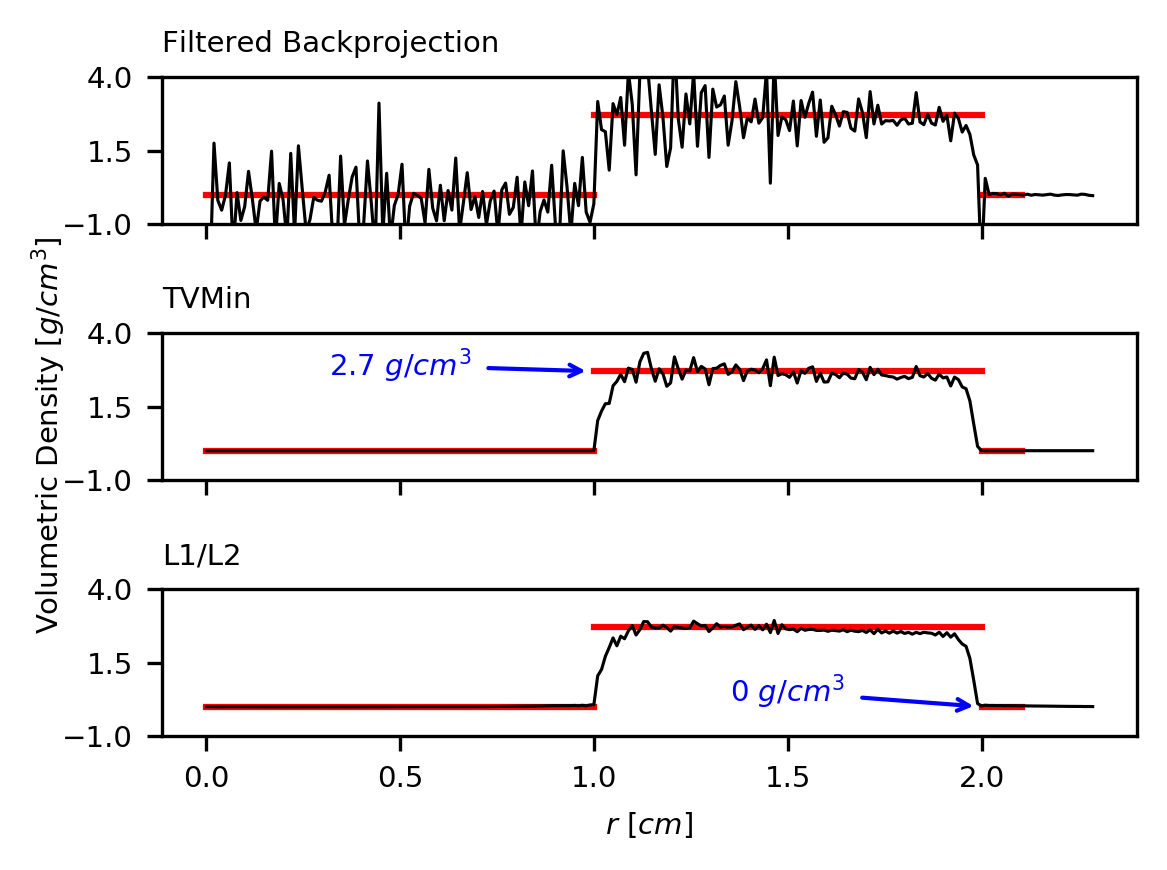}
    \caption{A comparison of horizontal samples of reconstructed volumetric density information using filtered backprojection (top), TVmin (center) and $L_1/L_2$ (bottom). All three plots correspond to $y=0$ plane, intersecting the origin. The red lines represent an estimated ground truth, given the density of pure Aluminum is 2.7$g/cm^3$, and assuming the remainder of the reconstruction is empty space.}
    \label{fig:lineoutcompare}
\end{figure}

In Figure \ref{fig:compare3}, we highlight a region of the reconstructed volume containing both Aluminum and empty regions. Qualitatively, the filtered back-projection solution recovers strong features of the target object, but endures the low-$r$ distortions and noise amplifications common in single-view tomography. Algorithm \ref{alg:boxtv} proved to be effective at preserving a strong characterization of the boundary and local statistics, but Algorithm \ref{alg:box} produces a more uniform reconstruction. 
We continue this assessment in Figure \ref{fig:lineoutcompare} by selecting a horizontal row of pixels from the reconstructed images corresponding to the central optical axis. The red lines underlayed upon the plots corresponds to an estimate of a ground truth garnered from the machining specifications of the calibration target. Similar to Figure \ref{fig:compare3} we see that, as parameterized, the $L_1/L_2$ method is comparatively more effective at suppressing noise and other artifacts. 

\section{Conclusions}

We presented a numerical optimization approach to solve single view tomographic reconstruction problems that utilize a constrained $L_1/L_2$ regularization. We considered both discrete cone and parallel beam formulations. Using the alternating direction of multipliers, we found that a well-parameterized $L_1/L_2$ scheme can return high-quality reconstructions, even in the presence of substantial noise. We provided numerical verification of these results in both synthetic, idealized settings and X-ray radiography.

% BibTeX users please use one of
%\bibliographystyle{spbasic}      % basic style, author-year citations
%\bibliographystyle{spmpsci}      % mathematics and physical sciences

\section*{Acknowledgements}
This manuscript has been authored in part by Mission Support and Test Services, LLC, under Contract No. DE-NA0003624 with the U.S. Department of Energy, National Nuclear Security Administration (DOE-NNSA), NA-10 Office of Defense Programs, and supported by the Site-Directed Research and Development Program. The United States Government retains and the publisher, by accepting the article for publication, acknowledges that the United States Government retains a non-exclusive, paid-up, irrevocable, world-wide license to publish or reproduce the published content of this manuscript, or allow others to do so, for United States Government purposes. The U.S. Department of Energy will provide public access to these results of federally sponsored research in accordance with the DOE Public Access Plan (http://energy.gov/downloads/doe-public-access-plan). The views expressed in the article do not necessarily represent the views of the U.S. Department of Energy or the United States Government. DOE/NV/036240--1104.

\bibliographystyle{siamplain}
\bibliography{references}
\end{document}